\documentclass{amsproc}

\usepackage{amssymb}

\newtheorem{theorem}{Theorem}[section]

\newtheorem{corollary}[theorem]{Corollary}
\theoremstyle{definition}

\newtheorem{assumption}[theorem]{Assumption}
\newtheorem{question}[theorem]{Question}

\newtheorem{conjecture}[theorem]{Conjecture}
\theoremstyle{remark}

\numberwithin{equation}{section}

\newcommand{\abs}[1]{\lvert#1\rvert}

\newcommand{\R}{{\mathbb R}}
\newcommand{\Z}{{\mathbb Z}}
\newcommand{\FF}{{\mathcal F}}
\newcommand{\MM}{{\mathcal M}}

\newcommand{\Int}{{{\rm Int}}}
\newcommand{\Cl}{{{\rm Cl}}}

\newcommand{\AAA}{{\mathcal A}}

\newcommand{\Fix}{{\rm Fix}}

\newcommand{\BB}{{\mathcal B}}

\newcommand{\HH}{{\mathbb H}}

\begin{document}
\title[horocycle flows for foliations]
{Remarks on the horocycle flows for foliations by hyperbolic surfaces}

{}
\author{Shigenori Matsumoto}
\address{Department of Mathematics, College of
Science and Technology, Nihon University, 1-8-14 Kanda, Surugadai,
Chiyoda-ku, Tokyo, 101-8308 Japan
}
\email{matsumo@math.cst.nihon-u.ac.jp
}
\thanks{The author is partially supported by Grant-in-Aid for
Scientific Research (C) No.\ 25400096.}


\date{\today }
\begin{abstract}
We show that the horocycle flow associated with a foliation on a compact manifold by
 hyperbolic surfaces is minimal under certain conditions.

\end{abstract}

\maketitle

\section{Introduction}
In the 1936 paper \cite{H}, G. A. Hedlund showed the minimality
of the horocycle flows associated to  closed oriented hyperbolic
surfaces. In \cite{MV}, the authors
consider a problem of generalizing this fact to compact laminations by
hyperbolic surfaces. 

Throughout this paper, we work under the following assumption.

\begin{assumption}
 $M$ is a closed smooth\footnote{In
this paper, smooth means $C^\infty$.} manifold, and $\FF$ a 
codimension $q$ minimal
foliation on $M$ by hyperbolic
surfaces \end{assumption}

The latter  condition means that $\FF$ is a 2-dimensional
smooth foliation  
equipped with a continuous leafwise metric of curvature $-1$
and that all the leaves of $\FF$ are dense in  $M$.
 Let $\Pi:\hat M\to M$ be the unit tangent bundle of the foliation $\FF$.
The total space $\hat M$ admits a locally free action  of $PSL(2,\R)$.
Denote the orbit
foliation by $\AAA$. See Section 2 of \cite{MV} for more detail.
(The foliation $\AAA$ is denoted by $T^1\FF$ in \cite{MV}.)
In other words,  we have the geodesic flow $g^t$,
the stable horocycle flow $h_+^t$ and the unstable horocycle flow
$h_-^t$. They preserve leaves of $\AAA$ and satisfy
\begin{equation}
\label{e0}
g^t\circ h_+^s\circ g^{-t}=h_+^{se^{-t}},\ \ \ \
g^t\circ h_-^s\circ g^{-t}=h_-^{se^{t}}.
\end{equation}
This says that the flow $g^t$ is uniformly 
hyperbolic along the leaves of $\AAA$.
The flows $g^t$ and $h^t_\pm$ jointly define an action
of a closed subgroup $B_{\pm}$ of $PSL(2,\R)$, whose orbit foliation is 
denoted by
$\BB_\pm$. They are subfoliations of $\AAA$ transverse to each other
in a leaf of $\AAA$ and the intersection is the orbit foliation of
$g^t$. The group $B_\pm$ is isomorphic to the group of the orientation
preserving affine transformations on the real line.
There is an involution $J:\hat M\to\hat M$ sending a leafwise unit
tangent vector $\zeta\in\hat M$ to $-\zeta$. The involution $J$
maps the flow $g^t$ to $g^{-t}$; $Jg^tJ=g^{-t}$, and the flow  $h_\pm^t$ to
$h_\mp^{-t}$; $Jh_\pm^tJ=h_\mp^{-t}$. Therefore the minimality of $\BB_+$ (resp.\ $h_+^t$) is
equivalent to the minimality of $\BB_-$ (resp.\ $h_-^t$).
Notice also that the minimality of $\FF$ immediately implies the
minimality of $\AAA$.
On the other hand,  there are examples of minimal $\FF$
for which the foliations
$\BB_+$ are not minimal \cite{MV}. The purpose of this paper is to study
the following question which arizes naturally from the result of Hedlund \cite{H}.

\begin{question}\label{question}
 Does the minimality of  the foliation
 $\BB_+$ 
imply the minimality
of the flow $h^t_+$?
\end{question}

So far no counter-examples are known.
There are some positive partial answers connected with this question.
(Below $(M,\FF)$ is assumed to satisfy Assumption 1.1.)

\medskip\noindent
(A) If $\FF$ is a homogeneous Lie foliation, then the flow $h_+^t$
is minimal \cite{ACD}.

\medskip\noindent
(B) If a leaf of $\FF$ admits a simple closed geodesic with
trivial holonomy, then the flow $h_+^t$ is minimal \cite{ADMV}. 

\medskip\noindent
(C) If the foliation $\BB_+$ is minimal and if there are nonplanar
leaf in $\FF$, then the flow $h_+^t$ admits a dense orbit \cite{MV}.

\medskip

Assume there is a simple closed oriented geodesic $c$ in some leaf of $\FF$.
Let $ D'\subset D$ be smooth closed $q$-disks in $M$ transverse
to $\FF$ such that $ D'\cap  c = D\cap  c=\{z_0\}$. (Recall that $\FF$
is of codimension $q$.) Let 
$ f: D'\to  D$ be the holonomy map of $\FF$ along the curve 
$c$. Thus $ f(z_0)=z_0$.
Let us consider the following condition.
See 
Figure 1.

\begin{figure}[h]
{\unitlength 0.1in%
\begin{picture}( 31.6200, 21.7400)( 18.0800,-26.3100)%
%
\special{pn 4}%
\special{sh 1}%
\special{ar 3370 1551 16 16 0  6.28318530717959E+0000}%
\special{sh 1}%
\special{ar 3370 1551 16 16 0  6.28318530717959E+0000}%
%
\special{pn 20}%
\special{ar 3360 1544 1552 1087  0.0000000  6.2831853}%
%
\special{pn 20}%
\special{ar 3370 1551 1330 931  0.2194833  0.2121438}%
%
\special{pn 8}%
\special{ar 2840 1551 504 353  0.0000000  6.2831853}%
%
\special{pn 8}%
\special{pa 3350 1551}%
\special{pa 3345 1516}%
\special{pa 3341 1482}%
\special{pa 3335 1448}%
\special{pa 3329 1415}%
\special{pa 3321 1383}%
\special{pa 3312 1352}%
\special{pa 3300 1323}%
\special{pa 3287 1296}%
\special{pa 3271 1271}%
\special{pa 3252 1249}%
\special{pa 3231 1229}%
\special{pa 3206 1211}%
\special{pa 3180 1196}%
\special{pa 3151 1183}%
\special{pa 3120 1172}%
\special{pa 3088 1162}%
\special{pa 3055 1154}%
\special{pa 2985 1142}%
\special{pa 2913 1134}%
\special{pa 2841 1130}%
\special{pa 2770 1130}%
\special{pa 2735 1132}%
\special{pa 2667 1140}%
\special{pa 2635 1147}%
\special{pa 2603 1155}%
\special{pa 2573 1165}%
\special{pa 2544 1176}%
\special{pa 2516 1190}%
\special{pa 2490 1205}%
\special{pa 2466 1223}%
\special{pa 2444 1243}%
\special{pa 2423 1265}%
\special{pa 2405 1289}%
\special{pa 2389 1316}%
\special{pa 2375 1344}%
\special{pa 2362 1374}%
\special{pa 2351 1405}%
\special{pa 2342 1437}%
\special{pa 2334 1470}%
\special{pa 2327 1504}%
\special{pa 2322 1539}%
\special{pa 2318 1573}%
\special{pa 2316 1607}%
\special{pa 2316 1641}%
\special{pa 2318 1674}%
\special{pa 2323 1706}%
\special{pa 2330 1737}%
\special{pa 2340 1766}%
\special{pa 2354 1793}%
\special{pa 2370 1817}%
\special{pa 2390 1839}%
\special{pa 2413 1859}%
\special{pa 2438 1878}%
\special{pa 2465 1894}%
\special{pa 2495 1908}%
\special{pa 2527 1921}%
\special{pa 2560 1933}%
\special{pa 2594 1943}%
\special{pa 2629 1952}%
\special{pa 2665 1961}%
\special{pa 2701 1968}%
\special{pa 2738 1974}%
\special{pa 2775 1979}%
\special{pa 2812 1983}%
\special{pa 2849 1986}%
\special{pa 2886 1987}%
\special{pa 2922 1987}%
\special{pa 2958 1985}%
\special{pa 2992 1982}%
\special{pa 3026 1977}%
\special{pa 3059 1971}%
\special{pa 3091 1963}%
\special{pa 3121 1953}%
\special{pa 3150 1941}%
\special{pa 3177 1927}%
\special{pa 3202 1912}%
\special{pa 3225 1894}%
\special{pa 3246 1874}%
\special{pa 3265 1852}%
\special{pa 3281 1828}%
\special{pa 3294 1802}%
\special{pa 3306 1774}%
\special{pa 3316 1744}%
\special{pa 3324 1713}%
\special{pa 3331 1681}%
\special{pa 3337 1648}%
\special{pa 3342 1614}%
\special{pa 3346 1580}%
\special{pa 3350 1551}%
\special{fp}%
%
\special{pn 4}%
\special{pa 3480 1075}%
\special{pa 3210 1180}%
\special{fp}%
\special{sh 1}%
\special{pa 3210 1180}%
\special{pa 3279 1174}%
\special{pa 3260 1161}%
\special{pa 3265 1137}%
\special{pa 3210 1180}%
\special{fp}%
\put(41.1000,-16.6300){\makebox(0,0)[lb]{$D'$}}%
%
\special{pn 4}%
\special{pa 4890 851}%
\special{pa 4730 921}%
\special{fp}%
\special{sh 1}%
\special{pa 4730 921}%
\special{pa 4799 913}%
\special{pa 4779 900}%
\special{pa 4783 876}%
\special{pa 4730 921}%
\special{fp}%
\put(49.7000,-8.7900){\makebox(0,0)[lb]{$D$}}%
\put(34.6000,-15.8600){\makebox(0,0)[lb]{$z_0$}}%
\put(26.7000,-15.9300){\makebox(0,0)[lb]{$U$}}%
\put(35.7000,-10.9600){\makebox(0,0)[lb]{$f(U)$}}%
\end{picture}}%
 \caption{}
\end{figure}
\medskip
\noindent
($*$) \
There exists an open subset $ U\subset  D'\subset D$ such that
$z_0\in\Cl( U)$, $ f( U)\supseteq  U$, and for any 
$z\in \Cl( U)$,  $d( f(z),z_0)\geq d(z,z_0)$.

\medskip

The main result of this paper is the following.

\begin{theorem} \label{main}
  If there is a leafwise simple closed geodesic which satisfies {\em($*$)},
then Question \ref{question} has a positive answer.
\end{theorem}

Notice that any leafwise simple closed geodesic with trivial holonomy
satisfies ($*$). 
So our result overlaps with the result (B) above.
Our proof, quite different from that of \cite{ADMV},
uses only the leafwise hyperbolicity (1.1) of the flow $g^t$.
More generally, any leafwise simple closed geodesic of
a Riemannian foliation satisfies ($*$). 
In fact, one can take as $U$ in ($*$) a transverse metric ball centered
at $z_0$.
On the other hand, it is shown
(Theorem 6, \cite{MV})
that if $\FF$ admits a holonomy invariant transverse measure, then
the foliation $\BB_+$ is minimal. Of course a Riemannian foliation
satisfies this condition. Thus we have the following corollary.

\begin{corollary}\label{riemannian}
 If $\FF$ is a 
 Riemannian foliation which admits a nonplanar leaf, then the flow 
$h_\pm^t$is minimal.
\end{corollary}

Unfortunately, this result, even combined with result (A), is not
sufficient to solve the following conjecture.

\begin{conjecture}
 For any Riemannian foliation $\FF$, the flow $h_\pm^t$ is minimal.
\end{conjecture}

Next consider the case where $\FF$ is of codimension one, that is,
$\dim(M)=3$.
We have a satisfactory answer in this case.
As will be shown in Section 4, a codimension one
foliation $\FF$ by hyperbodlic surfaces
 must have a nonplanar leaf. Let $ c$ 
be a leafwise simple closed geodesic. Consider the holonomy map $ f$
along $ c$ on a transverse open interval $I$ which intersects $c$ at a
point $z_1$. If the fixed point set $\Fix(f)$ has a nonempty
interior, then a point $z_0$ from the interior satisfies ($*$).
Otherwise an endpoint $z_0$ of a connected component $I\setminus\Fix(f)$
satisfies ($*$).
That is, if there is a nonplanar leaf, then ($*$) is always satisfied.
See Figure 2.
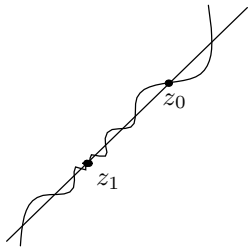
\begin{figure}[h]
{\unitlength 0.1in%
\begin{picture}( 12.6000, 12.6000)(  7.0000,-19.0000)%
%
\special{pn 8}%
\special{pa 700 1876}%
\special{pa 1960 652}%
\special{fp}%
%
\special{pn 8}%
\special{pa 1756 640}%
\special{pa 1766 720}%
\special{pa 1771 759}%
\special{pa 1773 798}%
\special{pa 1774 834}%
\special{pa 1773 869}%
\special{pa 1769 901}%
\special{pa 1761 931}%
\special{pa 1750 957}%
\special{pa 1735 980}%
\special{pa 1715 1000}%
\special{pa 1690 1015}%
\special{pa 1660 1026}%
\special{pa 1626 1034}%
\special{pa 1590 1040}%
\special{pa 1553 1045}%
\special{pa 1516 1051}%
\special{pa 1481 1058}%
\special{pa 1448 1068}%
\special{pa 1420 1081}%
\special{pa 1398 1100}%
\special{pa 1382 1124}%
\special{pa 1375 1156}%
\special{pa 1373 1228}%
\special{pa 1367 1259}%
\special{pa 1350 1278}%
\special{pa 1321 1285}%
\special{pa 1286 1287}%
\special{pa 1252 1291}%
\special{pa 1227 1307}%
\special{pa 1216 1337}%
\special{pa 1218 1375}%
\special{pa 1225 1410}%
\special{pa 1217 1430}%
\special{pa 1186 1430}%
\special{pa 1144 1423}%
\special{pa 1137 1440}%
\special{pa 1135 1467}%
\special{pa 1096 1468}%
\special{pa 1108 1468}%
\special{fp}%
%
\special{pn 8}%
\special{pa 772 1900}%
\special{pa 775 1860}%
\special{pa 779 1820}%
\special{pa 785 1782}%
\special{pa 792 1747}%
\special{pa 803 1715}%
\special{pa 817 1687}%
\special{pa 836 1664}%
\special{pa 860 1647}%
\special{pa 890 1636}%
\special{pa 926 1632}%
\special{pa 964 1631}%
\special{pa 1001 1628}%
\special{pa 1031 1618}%
\special{pa 1050 1595}%
\special{pa 1055 1557}%
\special{pa 1053 1513}%
\special{pa 1061 1487}%
\special{pa 1092 1496}%
\special{pa 1114 1502}%
\special{pa 1114 1473}%
\special{pa 1132 1468}%
\special{fp}%
%
\special{pn 4}%
\special{sh 1}%
\special{ar 1132 1468 16 16 0  6.28318530717959E+0000}%
\special{sh 1}%
\special{ar 1120 1468 16 16 0  6.28318530717959E+0000}%
\special{sh 1}%
\special{ar 1120 1468 16 16 0  6.28318530717959E+0000}%
%
\special{pn 4}%
\special{sh 1}%
\special{ar 1552 1048 16 16 0  6.28318530717959E+0000}%
\special{sh 1}%
\special{ar 1546 1048 16 16 0  6.28318530717959E+0000}%
\put(11.6800,-16.0000){\makebox(0,0)[lb]{$z_1$}}%
\put(15.1600,-11.8000){\makebox(0,0)[lb]{$z_0$}}%
\end{picture}}%

\caption{\small Although $z_1$ does not satisfy ($*$), some point
 $z_0\in I$  satisfies ($*$).}
\end{figure}
\\
Therefore we have:

\begin{corollary}
 \label{codim1} If $\FF$ is of codimension one, then Question
 \ref{question} has an affirmative answer.
\end{corollary}

Notice that there are many examples of the codimension one foliations $\FF$
for which $\BB_\pm$ are not minimal \cite{MV}.

\section{Proof of Theorem \ref{main}}

Let $ c$ be a leafwise simple closed geodesic in $M$ which satisfies
($*$). Associated to $ c$, there is a periodic orbit $\hat c$
of the leafwise
geodesic flow $g^t$ on $\hat M$, defined by $\hat c=( c, c')$, where $c'$ stands for
the derivative. 
Our overall strategy is to show that if $\BB_+$ is minimal, any minimal
set of $h_+^t$ intersects $\hat c$. For that purpose, we begin with
constructing a coordinate neighbourhood in a transversal of $\hat c$ for which 
the dynamics of the first return map of 
$g^t$ is well described. 

For a point $z_0$ in ($*$), let 
$\zeta_0\in\hat M$ be the tangent
vector of $c$ at $z_0$. Thus $\zeta_0\in\hat c$ and $\Pi(\zeta_0)=z_0$.
Let $\hat E$ be a smooth closed $(q+2)$-disk in the $(q+3)$-dimensional
manifold $\hat M$ transverse to
$g^t$ centered at $\zeta_0$.
Let $\hat D$ be a $q$-disk centered at $\zeta_0$ contained in $\Int(\hat
E)$
and transverse to the foliation $\AAA$.
If $\hat D$ is small enough, the projection $\Pi$ yields a diffeomorphism
from $\hat D$ to its image $\Pi(\hat D)$, which is a $q$-disk in $M$
transverse to $\FF$.
Here notice that the open set $U$ of ($*$) can be chosen arbitrarily
small (close
to $z_0$). In fact, if we replace $U$ by its intersection with
the metric disk  centered at $z_0$ of small radius in the transverse
$q$-disk, condition ($*$) is still satisfied. 
This shows that the $q$-disks $D$ and $D'$ in ($*$) can also
be chosen arbitrarily small. 
Therfore we may assume that
$D=\Pi(\hat D)$ is the disk in ($*$).

The disk $\hat E$, being 
transverse to the flow $g^t$, is transverse to the foliation $\BB_\pm$
and $\AAA$.
Let $\beta_{\pm}$ and $\alpha$ be the restriction of the foliation
$\BB_\pm$  and $\AAA$ to $\hat E$. See Figure 3.
\begin{figure}[h]
{\unitlength 0.1in%
\begin{picture}( 15.7800, 13.9600)(  4.7000,-20.9600)%
%
\special{pn 8}%
\special{pa 764 700}%
\special{pa 2038 700}%
\special{fp}%
\special{pa 774 705}%
\special{pa 470 1116}%
\special{fp}%
%
\special{pn 8}%
\special{pa 480 1112}%
\special{pa 1822 1112}%
\special{fp}%
%
\special{pn 8}%
\special{pa 2033 705}%
\special{pa 1827 1121}%
\special{fp}%
%
\special{pn 8}%
\special{pa 475 1112}%
\special{pa 475 2092}%
\special{fp}%
%
\special{pn 8}%
\special{pa 480 1758}%
\special{pa 774 1327}%
\special{fp}%
%
\special{pn 8}%
\special{pa 2038 705}%
\special{pa 2038 1685}%
\special{fp}%
%
\special{pn 8}%
\special{pa 475 2092}%
\special{pa 1847 2092}%
\special{fp}%
%
\special{pn 8}%
\special{pa 1842 2092}%
\special{pa 2038 1690}%
\special{fp}%
%
\special{pn 8}%
\special{pa 774 1327}%
\special{pa 2048 1327}%
\special{fp}%
%
\special{pn 8}%
\special{pa 480 1758}%
\special{pa 475 1758}%
\special{fp}%
%
\special{pn 8}%
\special{pa 480 1758}%
\special{pa 1852 1758}%
\special{fp}%
%
\special{pn 8}%
\special{pa 2048 1327}%
\special{pa 1832 1758}%
\special{fp}%
\special{pa 1832 1758}%
\special{pa 1832 1758}%
\special{fp}%
%
\special{pn 4}%
\special{pa 715 1410}%
\special{pa 1989 1410}%
\special{fp}%
%
\special{pn 4}%
\special{pa 666 1494}%
\special{pa 1940 1494}%
\special{fp}%
%
\special{pn 4}%
\special{pa 607 1572}%
\special{pa 1881 1572}%
\special{fp}%
%
\special{pn 4}%
\special{pa 543 1665}%
\special{pa 1817 1665}%
\special{fp}%
%
\special{pn 4}%
\special{pa 1989 1406}%
\special{pa 2009 1406}%
\special{fp}%
%
\special{pn 4}%
\special{pa 1930 1494}%
\special{pa 1969 1494}%
\special{fp}%
%
\special{pn 4}%
\special{pa 1881 1572}%
\special{pa 1930 1572}%
\special{fp}%
%
\special{pn 4}%
\special{pa 1822 1665}%
\special{pa 1876 1665}%
\special{fp}%
%
\special{pn 4}%
\special{pa 916 1332}%
\special{pa 656 1758}%
\special{fp}%
%
\special{pn 4}%
\special{pa 1058 1327}%
\special{pa 823 1758}%
\special{fp}%
%
\special{pn 4}%
\special{pa 1195 1337}%
\special{pa 989 1753}%
\special{fp}%
%
\special{pn 4}%
\special{pa 1940 1332}%
\special{pa 1734 1763}%
\special{fp}%
%
\special{pn 4}%
\special{pa 1813 1327}%
\special{pa 1626 1758}%
\special{fp}%
%
\special{pn 4}%
\special{sh 1}%
\special{ar 1274 1533 16 16 0  6.28318530717959E+0000}%
\special{sh 1}%
\special{ar 1278 1538 16 16 0  6.28318530717959E+0000}%
%
\special{pn 13}%
\special{pa 1278 1205}%
\special{pa 1278 1891}%
\special{fp}%
%
\special{pn 8}%
\special{pa 1832 1763}%
\special{pa 1832 2096}%
\special{fp}%
%
\special{pn 8}%
\special{pa 1832 1121}%
\special{pa 1832 1249}%
\special{fp}%
\put(19.6900,-20.2800){\makebox(0,0)[lb]{$\hat E$}}%
\put(12.2000,-20.5400){\makebox(0,0)[lb]{$\hat D$}}%
\put(7.6400,-18.7600){\makebox(0,0)[lb]{$\alpha(x)$}}%
\end{picture}}%
\caption{}
\end{figure}
\\ The $1$-dimensional foliations $\beta_+$ and $\beta_-$ are subfoliations
of 
the 2-dimensional foliation $\alpha$, transverse to each other in a leaf of
$\alpha$. Given a point $x\in \hat E$, let us denote by $\beta_\pm(x)$ and
 $\alpha(x)$
the leaves of the corresponding foliations which pass through $x$.
Given $\zeta\in \hat D$ and small $r>0$, let $\iota_{\pm,\zeta}:[-r,r]\to\beta_{\pm}(\zeta)$
be the isometric embedding such that $\iota_{\pm,\zeta}(0)=\zeta$. 
For $\xi,\eta\in[-r,r]$ and $\zeta\in\hat D$, let us denote
by $[\xi,\eta,\zeta]$ the unique point of the intersection of 
$\beta_+(\iota_{-,\zeta}(\xi))$ and $\beta_-(\iota_{+,\zeta}(\eta))$. 
See Figure 4.
\begin{figure}[h]
{\unitlength 0.1in%
\begin{picture}( 35.0000, 22.2400)( 13.9000,-27.2400)%
%
\special{pn 8}%
\special{pa 1390 1548}%
\special{pa 4780 1548}%
\special{fp}%
%
\special{pn 8}%
\special{pa 3010 500}%
\special{pa 3010 2724}%
\special{fp}%
%
\special{pn 20}%
\special{pa 1770 740}%
\special{pa 4330 740}%
\special{fp}%
%
\special{pn 20}%
\special{pa 1770 756}%
\special{pa 1770 2276}%
\special{fp}%
%
\special{pn 20}%
\special{pa 4330 2276}%
\special{pa 4330 740}%
\special{fp}%
%
\special{pn 20}%
\special{pa 1780 2284}%
\special{pa 4330 2284}%
\special{fp}%
%
\special{pn 4}%
\special{pa 3550 524}%
\special{pa 3550 2580}%
\special{fp}%
%
\special{pn 4}%
\special{pa 1640 1236}%
\special{pa 4650 1236}%
\special{fp}%
\put(48.9000,-16.1200){\makebox(0,0)[lb]{$\beta_-(\zeta)$}}%
\put(28.1000,-28.5200){\makebox(0,0)[lb]{$\beta_+(\zeta)$}}%
%
\special{pn 4}%
\special{sh 1}%
\special{ar 3010 1548 16 16 0  6.28318530717959E+0000}%
\special{sh 1}%
\special{ar 3010 1548 16 16 0  6.28318530717959E+0000}%
\put(22.5000,-17.1600){\makebox(0,0)[lb]{$\zeta=[0,0,\zeta]$}}%
\put(35.9000,-17.0000){\makebox(0,0)[lb]{$\xi$}}%
\put(29.1000,-11.6400){\makebox(0,0)[lb]{$\eta$}}%
\put(36.7000,-11.5600){\makebox(0,0)[lb]{$[\xi,\eta,\zeta]$}}%
\end{picture}}%
\caption{\small For each $\zeta\in\hat D$, the set $\{[\xi,\eta,\zeta]\mid \abs{\xi}\leq
 r,\abs{\eta}\leq r\}$ is a rectangle in $\alpha(\zeta)$.}
\end{figure}
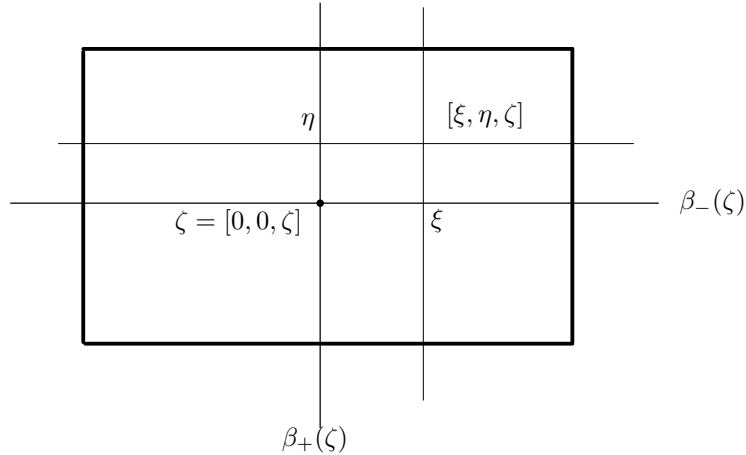

Recall the open subset $ U$ of $ D'$ in condition ($*$). For small $r>0$
(much smaller than the diameter of $U$), define
$$U_r=\{z\in U\mid d(z,z_0)< r\}.$$
 We have 
$f(U_r)\supseteq U_r$. Let 
$$\hat U_r=\Pi^{-1}(U_r)\cap \hat D\subset\hat M$$ 
and for $0<r'\leq r$, define
$$
\hat V_{r,r'}=\{\{[\xi,\eta,\zeta]\mid \abs{\xi}\leq r,\abs{\eta}\leq
r,\zeta\in \Cl(\hat U_{r'})\}.
$$
The first return map $F:\hat V_{r,r'}\to \hat E$  of the flow $g^t$
preserves the foliations $\alpha$ and $\beta_\pm$, and
therefore can be written as
$$
F[\xi,\eta,\zeta]=[\phi_\zeta(\xi),\psi_\zeta(\eta),f(\zeta)].$$
By some abuse, the conjugate of the map $f$ in ($*$)
by $\Pi\vert_{\hat D}:\hat D\to D$ is denoted here by $f$.
It satisfies
$f(\zeta_0)=\zeta_0$, $f(\hat U_{r'})\supseteq \hat U_{r'}$ and
$$ d(f(\zeta),f(\zeta_0))\geq d(\zeta,\zeta_0)),\
\forall
\zeta\in \Cl(\hat U_{r'}),$$
for an appropriate metric $d$.
By the hyperbolicity (\ref{e0}) of the leafwise geodesic flow $g^t$, there is
$\lambda\in(0,1)$ such that
\begin{equation}
 \label{e2}
d(\phi_\zeta(\xi),\phi_\zeta(\xi'))\geq\lambda^{-1}d(\xi,\xi'), \ \
\forall \xi,\xi'\in[-r,r],\ \forall \zeta\in\Cl(\hat U_{r'}),$$
$$d(\psi_\zeta(\eta),\psi_\zeta(\eta'))\leq\lambda d(\eta,\eta'), \ \
\forall \eta,\eta'\in[-r,r],\ \forall \zeta\in\Cl(\hat U_{r'}).
\end{equation}
On the other hand, since $\zeta_0=[0,0,\zeta_0]$ is a fixed point of $F$, we have
for small $r>0$ and even smaller $r'=r'(r)>0$, if $\zeta\in\Cl(\hat U_{r'})$,
$$\phi_{f^{-1}(\zeta)}(-r)<-r<r<\phi_{f^{-1}(\zeta)}(r)$$
and
\begin{equation}
 \label{e3}
-r<\psi_{f^{-1}(\zeta)}(-r)<\psi_{f^{-1}(\zeta)}(r)<r. 
\end{equation}
Therefore 
$$
F(\hat V_{r,r'})\cap \hat V_{r,r'}=\{[\xi,\eta,\zeta]\mid \abs{\xi}\leq r,\
\psi_{f^{-1}(\zeta)}(-r)
\leq\eta\leq
\psi_{f^{-1}(\zeta)}(r), \ \zeta\in \hat U_{r'}\}.$$
Replacing $\zeta\in\Cl(\hat U_{r'})$ by $f^{-1}(\zeta)\in\Cl(\hat U_{r'})$ in (2.2), we get
$$
-r<\psi_{f^{-2}(\zeta)}(-r)<\psi_{f^{-2}(\zeta)}(r)<r.
$$
Applying $\psi_{f^{-1}(\zeta)}$, we have
$$
\psi_{f^{-1}(\zeta)}(-r)<\psi_{f^{-1}(\zeta)}\psi_{f^{-2}(\zeta)}(-r)
<\psi_{f^{-1}(\zeta)}\psi_{f^{-2}(\zeta)}(r)<\psi_{f^{-1}(\zeta)}(r).
$$
This way, inductive use of (2.2) shows
$$
-r<\psi_{f^{-1}(\zeta)}(-r)<\psi_{f^{-1}(\zeta)}\psi_{f^{-2}(\zeta)}(-r)<\psi_{f^{-1}(\zeta)}\psi_{f^{-2}(\zeta)}\psi_{f^{-3}(\zeta)}(-r)<\cdots
$$$$\cdots<\psi_{f^{-1}(\zeta)}\psi_{f^{-2}(\zeta)}\psi_{f^{-3}(\zeta)}(r)
<\psi_{f^{-1}(\zeta)}\psi_{f^{-2}(\zeta)}(r)<\psi_{f^{-1}(\zeta)}(r)<r
$$
and by the hyperbolicity (\ref{e2})
$$ \lim_{n\to\infty}\psi_{f^{-1}(\zeta)}\cdots\psi_{f^{-n}(\zeta)}(-r)
=\lim_{n\to\infty}\psi_{f^{-1}(\zeta)}\cdots\psi_{f^{-n}(\zeta)}(r)=:\Psi(\zeta).
$$
See Figure 5.
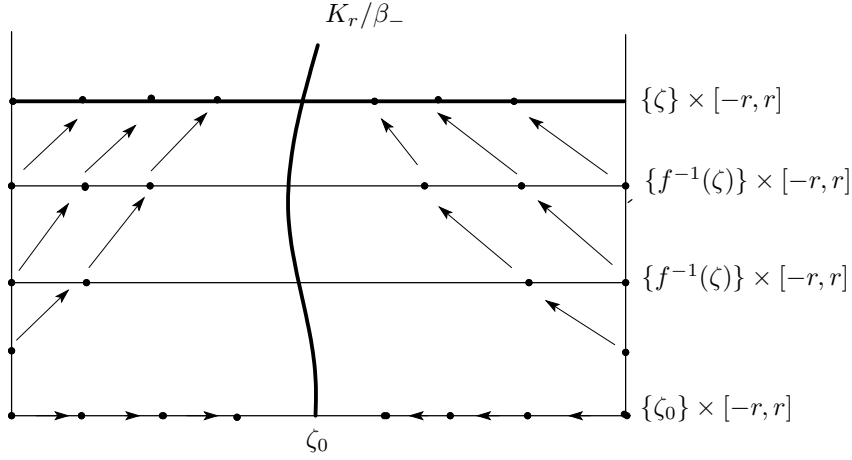
\begin{figure}
{\unitlength 0.1in%
\begin{picture}( 33.2400, 22.1800)( 15.9500,-27.5000)%
%
\special{pn 8}%
\special{pa 1195 2690}%
\special{pa 4429 2690}%
\special{fp}%
%
\special{pn 8}%
\special{pa 1208 2698}%
\special{pa 1208 678}%
\special{fp}%
%
\special{pn 8}%
\special{pa 4455 1558}%
\special{pa 4435 1573}%
\special{fp}%
%
\special{pn 8}%
\special{pa 4422 2690}%
\special{pa 4422 694}%
\special{fp}%
%
\special{pn 8}%
\special{pa 1215 1993}%
\special{pa 4422 1993}%
\special{fp}%
%
\special{pn 8}%
\special{pa 1208 1486}%
\special{pa 4422 1486}%
\special{fp}%
%
\special{pn 20}%
\special{pa 1208 1043}%
\special{pa 4415 1043}%
\special{fp}%
%
\special{pn 20}%
\special{pa 2811 749}%
\special{pa 2802 781}%
\special{pa 2793 812}%
\special{pa 2784 844}%
\special{pa 2775 875}%
\special{pa 2767 907}%
\special{pa 2758 939}%
\special{pa 2750 970}%
\special{pa 2741 1002}%
\special{pa 2733 1033}%
\special{pa 2726 1065}%
\special{pa 2718 1097}%
\special{pa 2711 1128}%
\special{pa 2704 1160}%
\special{pa 2697 1191}%
\special{pa 2685 1255}%
\special{pa 2680 1286}%
\special{pa 2675 1318}%
\special{pa 2671 1349}%
\special{pa 2667 1381}%
\special{pa 2663 1412}%
\special{pa 2660 1444}%
\special{pa 2658 1475}%
\special{pa 2657 1507}%
\special{pa 2655 1538}%
\special{pa 2655 1601}%
\special{pa 2656 1633}%
\special{pa 2658 1664}%
\special{pa 2661 1696}%
\special{pa 2664 1727}%
\special{pa 2668 1758}%
\special{pa 2672 1790}%
\special{pa 2677 1821}%
\special{pa 2683 1853}%
\special{pa 2689 1884}%
\special{pa 2695 1916}%
\special{pa 2701 1947}%
\special{pa 2708 1979}%
\special{pa 2722 2041}%
\special{pa 2729 2073}%
\special{pa 2736 2104}%
\special{pa 2750 2168}%
\special{pa 2756 2199}%
\special{pa 2763 2231}%
\special{pa 2769 2262}%
\special{pa 2775 2294}%
\special{pa 2780 2326}%
\special{pa 2785 2357}%
\special{pa 2793 2421}%
\special{pa 2796 2453}%
\special{pa 2798 2485}%
\special{pa 2800 2516}%
\special{pa 2801 2548}%
\special{pa 2801 2612}%
\special{pa 2799 2676}%
\special{pa 2799 2690}%
\special{fp}%
%
\special{pn 4}%
\special{sh 1}%
\special{ar 1215 1043 16 16 0  6.28318530717959E+0000}%
\special{sh 1}%
\special{ar 1215 1043 16 16 0  6.28318530717959E+0000}%
%
\special{pn 4}%
\special{sh 1}%
\special{ar 1580 1035 16 16 0  6.28318530717959E+0000}%
\special{sh 1}%
\special{ar 1580 1035 16 16 0  6.28318530717959E+0000}%
%
\special{pn 4}%
\special{sh 1}%
\special{ar 1939 1027 16 16 0  6.28318530717959E+0000}%
\special{sh 1}%
\special{ar 1939 1027 16 16 0  6.28318530717959E+0000}%
%
\special{pn 4}%
\special{sh 1}%
\special{ar 2285 1035 16 16 0  6.28318530717959E+0000}%
\special{sh 1}%
\special{ar 2285 1035 16 16 0  6.28318530717959E+0000}%
%
\special{pn 4}%
\special{sh 1}%
\special{ar 3108 1043 16 16 0  6.28318530717959E+0000}%
\special{sh 1}%
\special{ar 3108 1043 16 16 0  6.28318530717959E+0000}%
%
\special{pn 4}%
\special{sh 1}%
\special{ar 3441 1035 16 16 0  6.28318530717959E+0000}%
\special{sh 1}%
\special{ar 3441 1035 16 16 0  6.28318530717959E+0000}%
%
\special{pn 4}%
\special{sh 1}%
\special{ar 3838 1043 16 16 0  6.28318530717959E+0000}%
\special{sh 1}%
\special{ar 3838 1043 16 16 0  6.28318530717959E+0000}%
%
\special{pn 4}%
\special{sh 1}%
\special{ar 1208 1486 16 16 0  6.28318530717959E+0000}%
\special{sh 1}%
\special{ar 1208 1486 16 16 0  6.28318530717959E+0000}%
%
\special{pn 4}%
\special{sh 1}%
\special{ar 1593 1494 16 16 0  6.28318530717959E+0000}%
\special{sh 1}%
\special{ar 1593 1486 16 16 0  6.28318530717959E+0000}%
%
\special{pn 4}%
\special{sh 1}%
\special{ar 1933 1486 16 16 0  6.28318530717959E+0000}%
\special{sh 1}%
\special{ar 1933 1486 16 16 0  6.28318530717959E+0000}%
%
\special{pn 4}%
\special{sh 1}%
\special{ar 1208 1993 16 16 0  6.28318530717959E+0000}%
\special{sh 1}%
\special{ar 1208 1993 16 16 0  6.28318530717959E+0000}%
%
\special{pn 4}%
\special{sh 1}%
\special{ar 1600 1993 16 16 0  6.28318530717959E+0000}%
\special{sh 1}%
\special{ar 1600 1993 16 16 0  6.28318530717959E+0000}%
%
\special{pn 4}%
\special{sh 1}%
\special{ar 1208 2350 16 16 0  6.28318530717959E+0000}%
\special{sh 1}%
\special{ar 1208 2350 16 16 0  6.28318530717959E+0000}%
%
\special{pn 4}%
\special{sh 1}%
\special{ar 4422 1486 16 16 0  6.28318530717959E+0000}%
\special{sh 1}%
\special{ar 4422 1486 16 16 0  6.28318530717959E+0000}%
%
\special{pn 4}%
\special{sh 1}%
\special{ar 3877 1486 16 16 0  6.28318530717959E+0000}%
\special{sh 1}%
\special{ar 3877 1486 16 16 0  6.28318530717959E+0000}%
%
\special{pn 4}%
\special{sh 1}%
\special{ar 3370 1486 16 16 0  6.28318530717959E+0000}%
\special{sh 1}%
\special{ar 3370 1486 16 16 0  6.28318530717959E+0000}%
%
\special{pn 4}%
\special{sh 1}%
\special{ar 4422 1993 16 16 0  6.28318530717959E+0000}%
\special{sh 1}%
\special{ar 4422 1993 16 16 0  6.28318530717959E+0000}%
%
\special{pn 4}%
\special{sh 1}%
\special{ar 3916 1993 16 16 0  6.28318530717959E+0000}%
\special{sh 1}%
\special{ar 3916 1993 16 16 0  6.28318530717959E+0000}%
%
\special{pn 4}%
\special{sh 1}%
\special{ar 4422 2358 16 16 0  6.28318530717959E+0000}%
\special{sh 1}%
\special{ar 4422 2358 16 16 0  6.28318530717959E+0000}%
%
\special{pn 4}%
\special{sh 1}%
\special{ar 1208 2690 16 16 0  6.28318530717959E+0000}%
\special{sh 1}%
\special{ar 1208 2690 16 16 0  6.28318530717959E+0000}%
%
\special{pn 4}%
\special{sh 1}%
\special{ar 1574 2690 16 16 0  6.28318530717959E+0000}%
\special{sh 1}%
\special{ar 1574 2690 16 16 0  6.28318530717959E+0000}%
%
\special{pn 4}%
\special{sh 1}%
\special{ar 1997 2690 16 16 0  6.28318530717959E+0000}%
\special{sh 1}%
\special{ar 2003 2690 16 16 0  6.28318530717959E+0000}%
%
\special{pn 4}%
\special{sh 1}%
\special{ar 2388 2698 16 16 0  6.28318530717959E+0000}%
\special{sh 1}%
\special{ar 2388 2698 16 16 0  6.28318530717959E+0000}%
%
\special{pn 4}%
\special{sh 1}%
\special{ar 4429 2690 16 16 0  6.28318530717959E+0000}%
\special{sh 1}%
\special{ar 4415 2682 16 16 0  6.28318530717959E+0000}%
%
\special{pn 4}%
\special{sh 1}%
\special{ar 3504 2690 16 16 0  6.28318530717959E+0000}%
\special{sh 1}%
\special{ar 3504 2690 16 16 0  6.28318530717959E+0000}%
%
\special{pn 4}%
\special{sh 1}%
\special{ar 3159 2690 16 16 0  6.28318530717959E+0000}%
\special{sh 1}%
\special{ar 3171 2690 16 16 0  6.28318530717959E+0000}%
%
\special{pn 4}%
\special{sh 1}%
\special{ar 3909 2690 16 16 0  6.28318530717959E+0000}%
\special{sh 1}%
\special{ar 3909 2690 16 16 0  6.28318530717959E+0000}%
%
\special{pn 4}%
\special{pa 1337 2690}%
\special{pa 1504 2690}%
\special{fp}%
\special{sh 1}%
\special{pa 1504 2690}%
\special{pa 1437 2670}%
\special{pa 1451 2690}%
\special{pa 1437 2710}%
\special{pa 1504 2690}%
\special{fp}%
%
\special{pn 4}%
\special{pa 1804 2690}%
\special{pa 1875 2690}%
\special{fp}%
\special{sh 1}%
\special{pa 1875 2690}%
\special{pa 1808 2670}%
\special{pa 1822 2690}%
\special{pa 1808 2710}%
\special{pa 1875 2690}%
\special{fp}%
%
\special{pn 4}%
\special{pa 4237 2690}%
\special{pa 4127 2690}%
\special{fp}%
\special{sh 1}%
\special{pa 4127 2690}%
\special{pa 4194 2710}%
\special{pa 4180 2690}%
\special{pa 4194 2670}%
\special{pa 4127 2690}%
\special{fp}%
%
\special{pn 4}%
\special{pa 3761 2690}%
\special{pa 3646 2690}%
\special{fp}%
\special{sh 1}%
\special{pa 3646 2690}%
\special{pa 3713 2710}%
\special{pa 3699 2690}%
\special{pa 3713 2670}%
\special{pa 3646 2690}%
\special{fp}%
%
\special{pn 4}%
\special{pa 3389 2690}%
\special{pa 3293 2690}%
\special{fp}%
\special{sh 1}%
\special{pa 3293 2690}%
\special{pa 3360 2710}%
\special{pa 3346 2690}%
\special{pa 3360 2670}%
\special{pa 3293 2690}%
\special{fp}%
%
\special{pn 4}%
\special{pa 2125 2690}%
\special{pa 2285 2690}%
\special{fp}%
\special{sh 1}%
\special{pa 2285 2690}%
\special{pa 2218 2670}%
\special{pa 2232 2690}%
\special{pa 2218 2710}%
\special{pa 2285 2690}%
\special{fp}%
%
\special{pn 4}%
\special{pa 1593 1414}%
\special{pa 1856 1170}%
\special{fp}%
\special{sh 1}%
\special{pa 1856 1170}%
\special{pa 1794 1201}%
\special{pa 1817 1206}%
\special{pa 1821 1230}%
\special{pa 1856 1170}%
\special{fp}%
%
\special{pn 4}%
\special{pa 1939 1407}%
\special{pa 2228 1106}%
\special{fp}%
\special{sh 1}%
\special{pa 2228 1106}%
\special{pa 2167 1140}%
\special{pa 2191 1144}%
\special{pa 2196 1168}%
\special{pa 2228 1106}%
\special{fp}%
%
\special{pn 4}%
\special{pa 4326 1423}%
\special{pa 3921 1130}%
\special{fp}%
\special{sh 1}%
\special{pa 3921 1130}%
\special{pa 3963 1185}%
\special{pa 3964 1161}%
\special{pa 3987 1153}%
\special{pa 3921 1130}%
\special{fp}%
%
\special{pn 4}%
\special{pa 3852 1423}%
\special{pa 3441 1106}%
\special{fp}%
\special{sh 1}%
\special{pa 3441 1106}%
\special{pa 3482 1163}%
\special{pa 3483 1139}%
\special{pa 3506 1131}%
\special{pa 3441 1106}%
\special{fp}%
%
\special{pn 4}%
\special{pa 3332 1375}%
\special{pa 3145 1137}%
\special{fp}%
\special{sh 1}%
\special{pa 3145 1137}%
\special{pa 3170 1202}%
\special{pa 3178 1179}%
\special{pa 3202 1177}%
\special{pa 3145 1137}%
\special{fp}%
%
\special{pn 4}%
\special{pa 1247 1929}%
\special{pa 1504 1581}%
\special{fp}%
\special{sh 1}%
\special{pa 1504 1581}%
\special{pa 1448 1623}%
\special{pa 1472 1624}%
\special{pa 1480 1647}%
\special{pa 1504 1581}%
\special{fp}%
%
\special{pn 4}%
\special{pa 1613 1914}%
\special{pa 1901 1549}%
\special{fp}%
\special{sh 1}%
\special{pa 1901 1549}%
\special{pa 1844 1589}%
\special{pa 1868 1591}%
\special{pa 1875 1614}%
\special{pa 1901 1549}%
\special{fp}%
%
\special{pn 4}%
\special{pa 1234 2294}%
\special{pa 1504 2033}%
\special{fp}%
\special{sh 1}%
\special{pa 1504 2033}%
\special{pa 1442 2065}%
\special{pa 1466 2070}%
\special{pa 1470 2094}%
\special{pa 1504 2033}%
\special{fp}%
%
\special{pn 4}%
\special{pa 3825 1914}%
\special{pa 3447 1597}%
\special{fp}%
\special{sh 1}%
\special{pa 3447 1597}%
\special{pa 3485 1655}%
\special{pa 3488 1631}%
\special{pa 3511 1625}%
\special{pa 3447 1597}%
\special{fp}%
%
\special{pn 4}%
\special{pa 4351 1898}%
\special{pa 3967 1566}%
\special{fp}%
\special{sh 1}%
\special{pa 3967 1566}%
\special{pa 4004 1625}%
\special{pa 4007 1601}%
\special{pa 4031 1594}%
\special{pa 3967 1566}%
\special{fp}%
%
\special{pn 4}%
\special{pa 4365 2318}%
\special{pa 4006 2088}%
\special{fp}%
\special{sh 1}%
\special{pa 4006 2088}%
\special{pa 4051 2141}%
\special{pa 4051 2117}%
\special{pa 4073 2107}%
\special{pa 4006 2088}%
\special{fp}%
\put(44.9900,-11.0600){\makebox(0,0)[lb]{$\{\zeta\}\times[-r,r]$}}%
\put(45.0500,-20.4100){\makebox(0,0)[lb]{$\{f^{-1}(\zeta)\}\times[-r,r]$}}%
\put(45.1900,-15.3300){\makebox(0,0)[lb]{$\{f^{-1}(\zeta)\}\times[-r,r]$}}%
\put(28.5100,-6.6200){\makebox(0,0)[lb]{$K_r/\beta_-$}}%
\put(27.4800,-28.8000){\makebox(0,0)[lb]{$\zeta_0$}}%
\put(44.9900,-27.2100){\makebox(0,0)[lb]{$\{\zeta_0\}\times[-r,r]$}}%
%
\special{pn 4}%
\special{pa 1272 1400}%
\special{pa 1548 1140}%
\special{fp}%
\special{sh 1}%
\special{pa 1548 1140}%
\special{pa 1486 1171}%
\special{pa 1509 1177}%
\special{pa 1513 1200}%
\special{pa 1548 1140}%
\special{fp}%
\end{picture}}%
\caption{The map induced by $F$ on the quotient space $\hat V_r/\beta_-$.}
\end{figure}
Therefore for any $n>0$,
$$\bigcap_{0\leq i \leq n}F^i(\hat V_{r,r'})=\{[\xi,\eta,\zeta]\mid \abs{\xi}\leq
r,\
$$$$\psi_{f^{-1}(\zeta)}\cdots\psi_{f^{-n}(\zeta)}(-r)\leq\eta
\leq
\psi_{f^{-1}(\zeta)}\cdots\psi_{f^{-n}(\zeta)}(r), \ \zeta\in \Cl(\hat U_{r'})\},
$$
and
$$
K_{r,r'}:=\bigcap_{n\geq0}F^n(\hat V_{r,r'})=\{[\xi,\Psi(\zeta),\zeta]\mid
\abs{\xi}\leq r,
\zeta\in \Cl(\hat U_{r'})\}.$$
The subset $K_{r,r'}$ is closed, and by the closed graph theorem,
the function $\Psi:\Cl(\hat U_{r'})\to[-r,r]$ is continuous.
A crucial fact is that if $x\in K_{r,r'}$ and $n>0$, then $F^{-n}(x)\in \hat V_{r,r'}$.
(In fact, we have $F^{-n}(x)\in \hat K_{r,r'}$. But we do not use this.)

Let $\tau:\hat V_{r,r'}\to(0,\infty)$ be the first return time to $\hat E$ of the flow $g^t$
 and let
$$
\check V_{r,r'}=\{g^t(x)\mid x\in \hat V_{r,r'},t\in[0,\tau(x)]\}\ \mbox{ and }\
\check K_{r,r'}=\{g^t(x)\mid x\in K_{r,r'},t\in[0,\tau(x)]\}. $$
Both are compact sets, and we have:

\medskip\noindent
($*$$*$) \  If $x\in \check K_{r,r'}$ and $t>0$, then $g^{-t}(x)\in \check V_{r,r'}$.

\medskip
Now let us finish the proof of Theorem \ref{main}. 
We assume that the foliation $\BB_+$ is minimal. Let $\MM$ be any
minimal set of the flow $h_+^t$. Then we have
\begin{equation}
 \label{e23}
\bigcap_{t_0\in\R}\overline{\bigcup_{t\geq t_0}g^t(\MM)}=\hat 
M,
\end{equation}
since the LHS is $\BB_+$--invariant, closed and nonempty.
To show the $\BB_+$--invariance, we need to show that the LHS is
invariant both by $g^s$ and $h_+^s$. For the former, we have
$$
g^s(\bigcap_{t_0\in\R}\overline{\bigcup_{t\geq t_0}g^t(\MM)})
=\bigcap_{t_0\in\R}g^s\overline{\bigcup_{t\geq t_0}g^t(\MM)}
=\bigcap_{t_0\in\R}\overline{g^s(\bigcup_{t\geq t_0}g^t(\MM)})
$$$$=\bigcap_{t_0\in\R}\overline{\bigcup_{t\geq t_0}g^{s+t}(\MM)})
=\bigcap_{t_0\in\R}\overline{\bigcup_{t\geq t_0+s}g^t(\MM)}
=\bigcap_{t_0\in\R}\overline{\bigcup_{t\geq t_0}g^t(\MM)}.
$$
For the latter, notice that $g^t(\MM)$ is invariant by $h_+^s$ 
and therefore
$$
h_+^s(\bigcap_{t_0\in\R}\overline{\bigcup_{t\geq t_0}g^t(\MM)})
=\bigcap_{t_0\in\R}\overline{\bigcup_{t\geq t_0}h_+^s(g^t(\MM))}
=\bigcap_{t_0\in\R}\overline{\bigcup_{t\geq t_0}g^t(\MM)}.
$$

Now (\ref{e23}) implies in particular 
$\displaystyle \overline{\bigcup_{t\geq 0}g^t(\MM)}=\hat 
M,$
Since $\check V_{r,r'}$ has nonempty interior, 
we have $\displaystyle\bigcup_{t\geq 0}g^t(\MM)\cap\check V_{r,r'}\neq\emptyset$.
That is,
there is $x\in\MM$
and
$t\geq0$ such that $y=g^t(x)\in \check V_{r,r'}$.
Then an orbit segment of $h^t_+$ through $y$ intersects $\check K_{r,r'}$, say at a
point $y'$; $y'=h^s_+(y)\in \check K_{r,r'}$.
This is true for any $y\in \hat V_{r,r'}$ (Figure 5),
and for any $y\in \check V_{r,r'}$ by (1.1).
See Figures 6. 
\begin{figure}
{\unitlength 0.1in%
\begin{picture}( 23.5200, 12.6400)(  7.3000,-19.3400)%
%
\special{pn 8}%
\special{pa 730 670}%
\special{pa 730 1630}%
\special{dt 0.045}%
%
\special{pn 8}%
\special{pa 730 678}%
\special{pa 1482 1110}%
\special{dt 0.045}%
%
\special{pn 8}%
\special{pa 1482 1110}%
\special{pa 1482 1910}%
\special{dt 0.045}%
%
\special{pn 8}%
\special{pa 730 1574}%
\special{pa 1482 1934}%
\special{dt 0.045}%
%
\special{pn 8}%
\special{pa 1098 902}%
\special{pa 1098 1758}%
\special{dt 0.045}%
%
\special{pn 8}%
\special{pa 730 678}%
\special{pa 2650 678}%
\special{fp}%
%
\special{pn 8}%
\special{pa 1098 894}%
\special{pa 2858 894}%
\special{fp}%
%
\special{pn 8}%
\special{pa 1482 1118}%
\special{pa 3082 1118}%
\special{fp}%
%
\special{pn 8}%
\special{pa 1098 1750}%
\special{pa 2858 1750}%
\special{fp}%
%
\special{pn 8}%
\special{pa 1482 1926}%
\special{pa 3082 1926}%
\special{fp}%
%
\special{pn 8}%
\special{pa 2658 678}%
\special{pa 3066 1118}%
\special{dt 0.045}%
%
\special{pn 8}%
\special{pa 3066 1118}%
\special{pa 3066 1918}%
\special{dt 0.045}%
%
\special{pn 8}%
\special{pa 2858 894}%
\special{pa 2858 1750}%
\special{dt 0.045}%
%
\special{pn 4}%
\special{sh 1}%
\special{ar 1994 1374 16 16 0  6.28318530717959E+0000}%
\special{sh 1}%
\special{ar 1994 1374 16 16 0  6.28318530717959E+0000}%
%
\special{pn 4}%
\special{sh 1}%
\special{ar 2218 1534 16 16 0  6.28318530717959E+0000}%
\special{sh 1}%
\special{ar 2218 1534 16 16 0  6.28318530717959E+0000}%
%
\special{pn 8}%
\special{pa 1994 1374}%
\special{pa 2218 1534}%
\special{fp}%
\put(8.5000,-13.2600){\makebox(0,0)[lb]{$\check K_{r,r'}$}}%
\put(22.7400,-16.1400){\makebox(0,0)[lb]{$y$}}%
\put(18.5000,-15.1000){\makebox(0,0)[lb]{$y'$}}%
\end{picture}}%
\caption{}
\end{figure}
Let $x'=g^{-t}(y')$.
Then we have $$x'=g^{-t}(y')=g^{-t}h_+^s(y)=
g^{-t}h_+^sg^t(x)=h^{se^t}_+(x).$$
Hence $x'\in\MM$. On the other hand,
 since $y'\in\check K_{r,r'}$, we have
$x'\in \check V_{r,r'}$ by ($*$$*$). That is, $\MM\cap \check V_{r,r'}\neq\emptyset$.
Since $\bigcap_{r>0}\check V_{r,r'(r)}=\hat c$, we get by the finite intersection
property that

\medskip\noindent
($*$$*$$*$)\  $\MM\cap \hat c\neq\emptyset$.

\medskip 
The rest of the proof is routine.
For each $t\in\R$, we have either $\MM\cap g^t(\MM)=\emptyset$ or $\MM=g^t(\MM)$ since the both
sets are minimal sets of $h^t_+$. Let
$$
T=\{t\in\R\mid g^t(\MM)=\MM\}.$$
Then $T$ is a closed subgroup of $\R$. The statement ($*$$*$$*$) shows that
$T$ is nontrivial. If $T=\R$, then $\MM$ is $B_+$-invariant and we have
$\MM=\hat M$, as is required.
 Consider the remaining case where $T$ is isomorphic to $\Z$. In this
 case, the minimal set $\MM$ is a global cross section of $g^t$.
But this is impossible by Proposition 5 of \cite{MV}.
For the sake of completeness, let us give an easy proof for the case
where $M$ is a manifold.

The closed graph theorem shows that $\MM$ must be a tamely embedded topological
submanifold of codimension one. Thus the manifold $\hat M$ 
must be a bundle over $S^1$ and  admits a closed 1-form $\omega$
which takes positive value at $\displaystyle \frac{dg^t}{dt}(x)$ for any
$x\in \hat M$.
The closed geodesic $c$ which we started with and the closed
geodesic
with the reverse orientation correspond to
two periodic orbits $\hat c$ and $\hat c'$ of $g^t$ and we must have
$\displaystyle\int_{\hat c}\omega>0$ and $\displaystyle\int_{\hat c'}\omega>0$.
However $\hat c'$ is homotopic to $-\hat c$.
A contradiction.

\section{Codimension one foliations}
Here we shall prove the following.

\begin{theorem}
\label{last}
There is no foliation by heperbolic disks on any closed 3-manifold $M$.
\end{theorem}

{\sc Proof}. Assume on the contrary that there is a smooth foliation 
$\FF$ by hyperbolic disks on closed 3-manifold $M$.
H. Rosenberg \cite{R} showed that the 3-torus $T^3$ is the only
3-manifold which admits a smooth foliation by planes. So we have $M=T^3$.
According to W. Thurston \cite{T}, $\FF$ can be isotoped to be
transverse to a fibration $S^1\to T^3\to T^2$. 
Let $h:\Z^2=\pi_1(T^2)\to {\rm Diff}^\infty_+(S^1)$ be the holonomy 
homomorphism of the foliated bundle $\FF$. Then the associated $\Z^2$-action on $S^1$
 must be free,
 since all the leaves of $\FF$ are planar. 
Thus the action is topologically conjugate to an action by rigid rotations.
In particular, there is an $S^1$-action  on $S^1$ which
commutes with $h(\Z^2)$.
This implies that there is
an $\FF$-preserving topological $S^1$-action whose orbit foliation is the above
smooth fibration. Consider the covering space $S^1\times\R^2$ of $T^3$, where
the lifted foliation is $\{\{t\}\times \R^2\}$, and let
$\{\phi^t\}_{t\in S^1}$ be the lifted $S^1$-action. Each leaf
$\{t\}\times \R^2$ is equipped with a hyperbolic metric. Fix one leaf
$\{0\}\times\R^2$ which has a hyperbolic metric $g_0$ and replace
the metric of the other leaf $\{t\}\times\R^2$ by $(\phi^{-t})^*g_0$.
Then the new metric is $K$-quasiconformally equivalent to the old one
with fixed constant $K$. The quotient space of $S^1\times\R^2$ by
the $S^1$-action is identified
with the Poincar\'e upper half plane $\HH$. The covering transformation induces a $K$-quaiconformal
 action of $\Z^2$ on
$\HH$. Now a theorem of
D. Sullivan \cite{S} shows that such an action is topologically
conjugate to an action of a subgroup of $PSL(2,\R)$. Being a quotient action of
a covering transformation, this action must be cocompact, that is,
there is a compact subset of $\HH$ which intersects each orbit.
But this is impossible since the group is $\Z^2$, showing Theorem
\ref{last}. \qed

\end{document}